\documentclass[12pt]{amsart}

\newcommand{\GL}{\mathrm{GL}} %
\newcommand{\B}{\mathrm B} %
\renewcommand{\P}{\mathrm P} %
\newcommand{\U}{\mathrm U} %
\newcommand{\Part}{\mathbb P} %
\newcommand{\MultiPart}{\mathbb {MP}} %

\newcommand\BB{\mathbb B}
\newcommand\FF{\mathbb F}
\newcommand\ZZ{{\mathbb Z}}

\newcommand\CC{{\mathcal C}}
\newcommand\CE{{\mathcal E}}
\newcommand\CR{{\mathcal R}}
\newcommand\CQ{{\mathcal Q}}

\newcommand\sub{\subseteq}

\usepackage{amssymb}
\usepackage{fullpage}

\numberwithin{equation}{section}

\theoremstyle{plain}
\newtheorem{thm}[equation]{Theorem}
\newtheorem{prop}[equation]{Proposition}

\theoremstyle{definition}
\newtheorem{exmps}[equation]{Examples}

\theoremstyle{remark}

\title[Parabolic conjugacy in general linear groups]
{Parabolic conjugacy in general linear groups}
\author[S.~M.~Goodwin and G.~R\"ohrle]
{Simon M.~Goodwin and Gerhard R\"ohrle}

\address{School of Mathematics, University of Birmingham,
Birmingham, B15 2TT, UK} \email{goodwin@maths.bham.ac.uk}
\urladdr{http://web.mat.bham.ac.uk/S.M.Goodwin/}

\address{School of Mathematics, University of Southampton,
Southampton, SO17 1BJ, UK} \email{G.Roehrle@soton.ac.uk}
\urladdr{http://www.maths.soton.ac.uk/staff/Roehrle/}



\begin{document}

\begin{abstract}
Let $q$ be a power of a prime and $n$ a positive integer. Let $P(q)$
be a parabolic subgroup of the finite general linear group
$\GL_n(q)$. We show that the number of $P(q)$-conjugacy classes in
$\GL_n(q)$ is, as a function of $q$, a polynomial in $q$ with
integer coefficients. This answers a question of J.~Alperin in
\cite{A}.
\end{abstract}

\maketitle

\section{Introduction}
\label{sect:intro}

Let $\GL_n(q)$ be the general linear group of nonsingular $n \times
n$ matrices over the finite field $\FF_q$ and let $\U_n(q)$ be the
subgroup of $\GL_n(q)$ consisting of upper unitriangular matrices. A
longstanding conjecture states that the number of conjugacy classes
of $\U_n(q)$ is, as a function of $q$, a polynomial in $q$ with
integer coefficients. This conjecture has been attributed to
G.~Higman cf.\ \cite{H}, and it has been verified by computer for $n
\le 13$ in the work of A.~Vera-L\'opez and J.~M.~Arregi, see
\cite{VA}. There has been further interest in this conjecture from
G.~R.~Robinson \cite{Rob} and J.~Thompson \cite{T}.

In \cite{A}, J.~Alperin showed that a related result is ``easily
established'', namely that the number of $\U_n(q)$-conjugacy classes
in all of $\GL_n(q)$ is a polynomial in $q$ with integer
coefficients.  This theorem can be viewed as evidence in support of
Higman's conjecture.  Alperin also considers the possibility of a
proof of Higman's conjecture by descent from the theorem proved in
\cite{A}, though he says that this seems very unlikely.

In addition, Alperin showed in \cite{A} that the number of
$\B_n(q)$-conjugacy classes in $\GL_n(q)$ is a polynomial in $q$,
where $\B_n(q)$ is the subgroup of upper triangular matrices in
$\GL_n(q)$.

\medskip

Let $d = (d_1,\dots,d_t) \in \ZZ_{\ge 1}^t$ satisfy $d_i < d_{i+1}$
and $d_t = n$; we call such $d$ an {\em $n$-dimension vector}. Let
$\P_{n,d}(q)$ be the parabolic subgroup of $\GL_n(q)$ that
stabilizes the standard flag $\{0 \} \sub \FF_q^{d_1} \sub
\FF_q^{d_2} \sub \dots \sub \FF_q^{d_t} = \FF_q^n$, and let
$\U_{n,d}(q)$ be the unipotent radical of $\P_{n,d}(q)$.  In
\cite{A} Alperin asks if the number of $\U_{n,d}(q)$-conjugacy
classes in $\GL_n(q)$ is a polynomial in $q$; and likewise for the
number of $\P_{n,d}(q)$-conjugacy classes in $\GL_n(q)$. In
\cite[Thm.\ 4.5]{GR}, the authors showed that this question for
$\U_{n,d}(q)$ has an affirmative answer. In this paper, we prove the
following theorem, which affirmatively answers Alperin's question
for $\P_{n,d}(q)$.

\begin{thm}
\label{T:main}
The number of $\P_{n,d}(q)$-conjugacy classes in $\GL_n(q)$ is, as a
function of $q$ for fixed $d$, a polynomial in $q$ with integer
coefficients.
\end{thm}

The special case of Theorem \ref{T:main} when
$\P_{n,d}(q) = \GL_n(q)$ is well known, of course.

In order to state a proposition related to Theorem \ref{T:main}, we
need to recall some standard terminology.  We let $K$ be the
algebraic closure of $\FF_q$ and view $\GL_n(q)$ as a subgroup of
$\GL_n(K)$ in the natural way.  Recall that two parabolic subgroups
of $\GL_n(K)$ are said to be {\em associated} if they have Levi
subgroups that are conjugate in $\GL_n(K)$.  We write $\P_{n,d}(K)$
for the parabolic subgroup of $\GL_n(K)$ such that $\P_{n,d}(K) \cap
\GL_n(q) = \P_{n,d}(q)$. Let $d = (d_1,\dots,d_t)$ and $d' =
(d'_1,\dots,d'_{t'})$ be $n$-dimension vectors. We recall that
$\P_{n,d}(K)$ and $\P_{n,d'}(K)$ are associated if and only if $t =
t'$ and there exists $\sigma \in \mathrm{Sym}(t)$ such that $d_i -
d_{i-1} = d'_{\sigma i} - d'_{\sigma i - 1}$ for all $i =
1,\dots,t$; by convention we set $d_0 = d'_0 = 0$.

By \cite[(4.15)]{GR}, we have the following proposition.  We
indicate how it is proved in the outline of the proof of Theorem
\ref{T:main} given below.

\begin{prop} \label{P:assoc}
Let $\P_{n,d}(K)$ and $\P_{n,d'}(K)$ be associated parabolic
subgroups of $\GL_n(K)$. Then the number of $\P_{n,d}(q)$-conjugacy
classes in $\GL_n(q)$ is equal to the number of
$\P_{n,d'}(q)$-conjugacy classes in $\GL_n(q)$.
\end{prop}

We note that the proof of the observation in Proposition
\ref{P:assoc} does not yield a bijection between the two sets of
orbits.  It would be interesting to know if a bijection can be
defined in a natural way.

\medskip

Below we give an outline of our proof of Theorem \ref{T:main}.
Before doing this, we simplify our notation. We write $G =
\GL_n(q)$, $B = \B_n(q)$, and, for $d$ as above, we write $P =
\P_{n,d}(q)$. For a subgroup $H$ of $G$, we write $k(H,G)$ for the
number of $H$-conjugacy classes in $G$. Although this notation does
not show a dependence on $q$, we want to allow $q$ to vary and for
$G$, $B$, $P$ to define groups for each $q$; so, for example, it
makes sense to say that $k(P,G)$ is a polynomial in $q$. We write
${\bf G} = \GL_n(K)$ and $\bf P$ for the parabolic subgroup of $\bf
G$ corresponding to $P$.

For $x \in G$, we define $f_P^G(x)$ to be the number of conjugates
of $P$ containing $x$, i.e.\ $f_P^G(x) = |\{{}^yP \mid y \in G, x
\in {}^yP\}|$.  A counting argument as in \cite{A} (see also
\cite[\S 4.1]{GR}), along with the fact that $P = N_G(P)$, yields
\begin{equation}
\label{e:kPG} k(P,G)  =  \sum_{x \in \CR} f_P^G(x),
\end{equation}
where $\CR = \CR(P,G)$ is a set of representatives of the conjugacy
classes of $G$ that intersect $P$.  We note that if the conjugacy
class of $x \in G$  misses $P$, then $f_P^G(x) = 0$. Therefore, it
does no harm in \eqref{e:kPG} to sum over a set of representatives
$\CR = \CR(G)$ of all conjugacy classes of $G$.

From the proof of \cite[Lem.\ 3.2]{GR}, one can observe that for $x
\in G$, $f_P^G(x)$ only depends on $P$ up to the \emph{association
class} of $\bf P$, i.e.\ if ${\bf P}$ and ${\bf Q}$ are associated
parabolic subgroups of ${\bf G}$, then $f_P^G(x) = f_Q^G(x)$ for all
$x \in G$. This is a consequence of the fact that the Harish-Chandra
induction functor $R_L^G$ is independent of the parabolic subgroup
that contains the Levi subgroup $L$.  This observation is used to
deduce \cite[(4.15)]{GR} and so Proposition \ref{P:assoc}.

In \cite{A}, Alperin shows that $k(B,G)$ is a polynomial in $q$,
using the formula \eqref{e:kPG} for the case $P=B$. The proof of
this depends on partitioning the set $\CR(B,G)$ into a finite union
$\CR(B,G) = \CR_1 \cup \cdots \cup \CR_r$ independent of $q$ (though
some $\CR_i$ may be empty for small $q$) such that $f_B^G(x) =
f_B^G(y)$ if $x,y \in \CR_i$; and $|\CR_i|$ is a polynomial in $q$.
An inductive counting argument is used to show that $f_B^G(x_i)$ is
given by a polynomial in $q$, for $x_i \in \CR_i$.

In this paper we give an analogous decomposition $\CR(G) = \CR_1
\cup \cdots \cup \CR_r$; this partition is based on Jordan normal
forms. Again, this decomposition does not depend on $q$ (though some
$\CR_i$ may be empty for small $q$) and we show that $|\CR_i|$ is a
polynomial in $q$.  Let $x \in \CR_i$, for some $i$, with Jordan
decomposition $x = su$, and let $H = C_G(s)$.  We show that
$f_P^G(x)$ can be expressed as a sum of terms of the form
$f_Q^H(u)$, where $Q$ is a parabolic subgroup of $H$ of the form
${}^yP \cap H$ for some $y \in G$.  If $x' = s'u' \in \CR_i$, then
we have $u' = u$, and so we have $f_P^G(x') = f_P^G(x)$. We can
appeal to \cite[Thm.\ 3.10]{GR} to deduce that each $f_Q^H(u)$ is a
polynomial in $q$ and therefore, that $f_P^G(x)$ is a polynomial in
$q$. The key point in the proof that $f_Q^H(u)$ is a polynomial in
$q$ is to show that it can be expressed in terms of Green functions;
in the present setting the results in \cite{Gr} show that these
Green functions are polynomials in $q$. We then have
\begin{equation}
\label{e:kPG2} k(P,G) = \sum_{i=1}^r |\CR_i| \: f_P^G(x_i),
\end{equation}
where $x_i \in \CR_i$.  Each summand on the right-hand side of
\eqref{e:kPG2} is a polynomial in $q$.  Hence, $k(P,G)$ is a
polynomial in $q$.

We are left to show that, as a polynomial in $q$, $k(P,G)$ has
integer coefficients. This is non-trivial: although the coefficients
of the polynomial $f_P^G(x)$ are integers (this follows from the
results in \cite[\S 4]{GR}), the coefficients of the polynomials
$|\CR_i|$ are not integers in general.  In order to show that
$k(P,G) \in \ZZ[q]$, we argue that the $P$-conjugacy classes in $G$
can be parameterized by the $\FF_q$-rational points of a family of
varieties defined over $\FF_q$.  Then we apply some standard
arguments.

\medskip

Let $U$ be the unipotent radical of $P$ and let $u \in G$ be
unipotent.  Using the theory of Green functions, it is proved in
\cite{GR} that $f_U^G(u)$ is a polynomial of $q$; also in the
appendix of {\em loc.\ cit.\ }an elementary counting argument is
used to give an alternative proof of this. It is possible to give an
elementary proof that $f_P^G(u)$ is a polynomial in $q$ for
unipotent $u$; this proof is similar to that in the appendix to
\cite{GR} and is rather technical, so we choose not to include it
here. Given such a proof one can avoid appealing to the theory of
Green functions in the proof of Theorem \ref{T:main}. For this one
needs to observe that for semisimple $s \in G$, the centralizer $H =
C_G(s)$ is isomorphic to a direct product of groups of the form
$\GL_{m}(q^l)$, where $m,l \in \ZZ_{\ge 1}$.  Then, for arbitrary $x
\in G$ with Jordan decomposition $x = su$, one can deduce that
$f_P^G(x)$ is a polynomial in $q$ using the expression for
$f_P^G(x)$ as a sum of terms of the form $f_Q^H(u)$.

In analogy to a comment made at the end of the appendix to
\cite{GR}, it is not possible to deduce Proposition \ref{P:assoc}
from an elementary proof of Theorem \ref{T:main} as described above.

\medskip

One can consider the more general situation where the general linear
group $\GL_n(q)$ is replaced by an arbitrary finite group of Lie
type $G$, and $P$ is a parabolic subgroup of $G$ with unipotent
radical $U$. The precise formulation of the analogous questions
regarding $k(U,G)$ and $k(P,G)$ being polynomials in $q$ with
integer coefficients is rather technical so we do not give it here;
this formulation requires an axiomatic setup as in \cite[\S
2.2]{GR}. However, we note that \cite[Thm.\ 4.5]{GR} says that
$k(U,G)$ is a polynomial in $q$ if $p$ is good for ${\bf G}$ and
${\bf G}$ has connected centre, where ${\bf G}$ is the connected
reductive algebraic group defined over $\FF_q$ so that $G$ is the
group of $\FF_q$-rational points of $\bf G$.  In case ${\bf G}$ has
disconnected centre, $k(U,G)$ is only given by polynomials up to
congruences on $q$. That is, in the language of G.~Higman \cite{H2},
$k(U,G)$ is PORC (Polynomial On Residue Classes); this is discussed
before \cite[Exmp.\ 4.10]{GR}. The question about $k(P,G)$ is more
difficult in general.  We believe that one should be able to
generalize the arguments in this paper to show that $k(P,G)$ is PORC
in general. As is mentioned in \cite[Rem.\ 4.12]{GR}, the centre of
a pseudo-Levi subgroup of ${\bf G}$ need not be connected even if
the centre of ${\bf G}$ is connected; therefore, in general one can
only hope to prove that $f_P^G(x)$ is PORC.

\smallskip

As a general reference for algebraic groups defined over finite
fields we refer the reader to the book by Digne and Michel
\cite{DM}.

\section{Notation}
\label{sect:notn}

We establish the notation to be used throughout this note.  We
continue to use the convention that the objects that we define
depend on the prime power $q$, but that this dependence is
suppressed in our notation.

\smallskip

We write $\FF_q$ for the finite field of $q$ elements. We denote the
algebraic closure of $\FF_q$ by $K$ and we consider all the finite
fields $\FF_{q^m}$ (for $m \in \ZZ_{\ge 1}$) as subfields of $K$.
The set of non-zero elements of $K$ is denoted by $K^\times$;
likewise $\FF_q^\times$ denotes the set of non-zero elements of
$\FF_q$. For $a \in K^\times$, the {\em degree of $a$ over $q$},
denoted $\deg(a) = \deg_q(a)$, is the minimal value of $m$ such that
$a \in \FF_{q^m}$. For $m \in \ZZ_{\ge 2}$ we define
$\FF_{q^m}^\sharp$ by
\[
\FF_{q^m}^\sharp = \FF_{q^m} \setminus \bigcup_{j|m} \FF_{q^j} = \{a
\in K \mid \deg(a) = m\};
\]
we define $\FF_q^\sharp = \FF_q^\times$.

We write $F$ for the Frobenius morphism on $K$ corresponding to $q$,
i.e.\ $F(a) = a^q$ for all $a \in K$. We let $K^\times/F$ denote the
set of $F$-orbits in $K^\times$; this set is in bijection with the
set of all monic irreducible polynomials in $\FF_q[X] \setminus
\{X\}$.  Given $a \in K$ we write $\bar a$ for the $F$-orbit of $a$
in $K$.  Note that the degree function is constant on $F$-orbits in
$K^\times$, so that for  given $\bar a \in K^\times/F$, the degree
$\deg(a)$ is well-defined. Also, we sometimes consider a sum or
product over $K^\times/F$ where the summands or factors are indexed
by representatives of the $F$-classes in $K^\times$; in such
situations each summand or factor only depends on the corresponding
element in $K^\times/F$.

Given a map $\gamma : K^\times/F \to S$, where $S$ is some
set, we write $\gamma_0 : K^\times \to S$ for the map defined by
$\gamma_0(a) = \gamma(\bar a)$.  For $m \in \ZZ_{\ge 1}$, we write
$\FF_{q^m}^\sharp/F$ for the set of $F$-orbits in
$\FF_{q^m}^\sharp$ and define
\begin{equation}
\label{e:Fsharp}
\phi(m) = |\FF_{q^m}^\sharp/F|.
\end{equation}
We observe that
\[
\phi(m) = \frac{1}{m} \sum_{j | m} \mu(j) q^{m/j},
\]
where $\mu$ is the classical M\"obius function, see for example
\cite[\S 1.13]{K}; in particular, $\phi(m)$ is a polynomial in $q$.

\smallskip

By a partition we mean a sequence of the form $\lambda =
(\lambda_1^{c_1},\dots,\lambda_l^{c_l})$, where $\lambda_i,c_i \in
\ZZ_{\ge 1}$ and $\lambda_i > \lambda_{i+1}$; we allow $\lambda$ to
be the empty partition, i.e.\ $l = 0$, $\lambda = ()$. Given a
partition $\lambda$, we let $|\lambda| = \sum_{i=1}^l c_i
\lambda_i$.  We write $\Part$ for the set of all partitions.

We fix a linear order $\prec$ on $\Part$, by setting $\lambda \prec
\lambda'$ if $|\lambda| < |\lambda'|$ and then ordering the
partitions $\lambda$ for fixed $|\lambda|$ lexicographically. By a
multi-partition we mean a sequence of the form $\mu =
(\mu_1^{b_1},\dots,\mu_m^{b_m})$, where $\mu_i \in \Part$, $b_i \in
\ZZ_{\ge 1}$ and $\mu_i \succ \mu_{i+1}$; we allow $\mu$ to be the
empty multi-partition. Given a multi-partition $\mu=
(\mu_1^{b_1},\dots,\mu_m^{b_m})$ we let $|\mu| =
\sum_{i=1}^m b_i|\mu_i|$. We write $\MultiPart$ for the set of all
multi-partitions.

\smallskip

The polynomial defined below is required to simplify notation in
Section \ref{sect:conj}. For a sequence $b=(b_1,\dots,b_m) \in
\ZZ_{\ge 1}^m$ we define the following polynomial in the
indeterminate $z$:
\begin{equation} \label{e:Delta}
\Delta(b,z) =
\binom{z}{b_1}\binom{z-b_1}{b_2}\binom{z-b_1-b_2}{b_3}\cdots
\binom{z-b_1-\dots -b_{m-1}}{b_m},
\end{equation}
where $\binom{z}{c} = \frac{z(z-1)\cdots(z-c+1)}{c!}$, for $c \in
\ZZ_{\ge 1}$.  We allow $\Delta$ to be defined for different values
of $m$. We note that the coefficients of $\Delta(b,z)$ are in
general not integers.

\smallskip

Let $n$ be a positive integer.  We write $G = \GL_n(q)$ and regard
it as a subgroup of ${\bf G} = \GL_n(K)$.  We write $F$ for the
standard Frobenius morphism on ${\bf G}$ and its natural module
$K^n$. Therefore, $G = {\bf G}^F$ is the group of fixed points of
$F$ in ${\bf G}$, and $\FF_q^n = (K^n)^F$.

For $g,x \in G$, we write ${}^g x = gxg^{-1}$; similarly for a
subgroup $H$ of $G$ we write ${}^g H = gHg^{-1}$. We write $C_G(x) =
\{g \in G \mid {}^g x = x\}$ for the centralizer of $x$ in $G$; the
centralizer of $x$ in ${\bf G}$ is denoted by $C_{{\bf G}}(x)$.

Let $m \in \ZZ_{\ge 1}$ and $a \in K$.  Then the $m \times m$ Jordan
matrix $J(a,m)$ is defined as usual.  Given a partition $\lambda
=(\lambda_1^{c_1},\dots,\lambda_l^{c_l})$, the matrix $J(a,\lambda)$
is defined as a direct sum of Jordan matrices:
\[
J(a,\lambda) = \bigoplus_{i=1}^l c_i J(a,\lambda_i).
\]
Finally, for $\bar a \in K^\times/F$ and $\lambda \in \Part$,
we define the matrix
\[
J(\bar a,\lambda) = \bigoplus_{i=0}^{\deg(a) - 1} J(F^i(a),\lambda).
\]
By choosing a basis of the form $\BB_0 \cup \BB_1 \cup \dots \cup
\BB_{\deg(a) -1}$ for $K^n$ (where $n = \deg(a) |\lambda|$) with
$|\BB_i| = |\lambda|$ and $F^i(\BB_0) = \BB_i$, the matrix $J(\bar
a,\lambda)$ is fixed by $F$ and so lies in $G$.

\section{The conjugacy classes of $\GL_n(q)$}
\label{sect:conj}

In this section we recall the parametrization of the conjugacy
classes of $G = \GL_n(q)$, see for example \cite[Ch.\ IV \S 2]{M}.
We use this parametrization to define the partition of the set of
conjugacy classes of $G$ mentioned in the introduction.

The conjugacy classes of $G$ are given by Jordan normal forms and
these are parameterized by  maps
\[
\gamma : K^\times/F \to \Part
\]
such that $\gamma(\bar a)$ is the empty partition for all but
finitely many $\bar a \in K^\times/F$ and
\[\sum_{a \in K^\times} |\gamma_0(a)| =
\sum_{\bar a \in K^\times/F} \deg(a) |\gamma(\bar a)| = n.
\]
We write
$\Gamma$ for the set of all such maps $\gamma$. Given $\gamma \in
\Gamma$, we can define a linear map $x(\gamma) \in G$ as follows: We
decompose $K^n$ as
\[
K^n = \bigoplus_{a \in K^\times} V_a,
\]
where $\dim V_a = |\gamma_0(a)| = |\gamma(\bar a)|$
and $F(V_a) = V_{F(a)}$ for all $a
\in K^\times$.  For $\bar a \in K^\times/F$, we write $V_{\bar a} =
\bigoplus_{i=0}^{\deg(a) -1} V_{F^i(a)}$.  With respect to an
(ordered) basis, denoted $\BB(\gamma)_{\bar a}$, of $V_{\bar a}$,
the action of $x(\gamma)$ on $V_{\bar a}$ is given by the matrix
$J(\bar a,\gamma(\bar a))$. The set $\{x(\gamma) \mid \gamma \in
\Gamma\}$ gives a complete set of representatives of the conjugacy
classes of $G$.

For $a \in K^\times$, we define $\BB(\gamma)_a =
\BB(\gamma)_{\bar a} \cap V_a$.  We write $\BB(\gamma)$ for the
basis of $K^n$ given by $\BB(\gamma) = \bigcup_{a \in K^\times}
\BB(\gamma)_a$.

\smallskip

Let $\gamma \in \Gamma$. We write the Jordan decomposition of
$x(\gamma)$ as $x(\gamma) = s(\gamma)u(\gamma)$. It is
straightforward to describe the action of $s(\gamma)$ and
$u(\gamma)$ on each $V_a$, for $a \in K^\times$.

The semisimple part $s(\gamma)$ acts on $V_a$ as multiplication by
$a$. Therefore, we see that the centralizer of $s(\gamma)$ in ${\bf
G}$ is
\[
C_{{\bf G}}(s(\gamma)) = \prod_{a \in K^\times} \GL(V_a) \cong
\prod_{\bar a \in K^\times/F} \GL_{|\gamma(\bar a)|}(K)^{\deg(a)}.
\]
In order to describe the centralizer of $s(\gamma)$ in $G$, we note
that $V_a$ is defined over $\FF_{q^{\deg(a)}}$, and
$V_a^{F^{\deg(a)}} \cong \FF_{q^{\deg(a)}}^{|\gamma_0(a)|}$.
Note that for $a, b \in K^\times$ in the same $F$-orbit, we have
$V_a^{F^{\deg(a)}} \cong V_b^{F^{\deg(b)}}$.
Therefore, as $F(V_a) = V_{F(a)}$, we see that the centralizer of
$s(\gamma)$ in $G$ is
\begin{equation} \label{e:cent}
C_G(s(\gamma)) \cong \prod_{\bar a \in K^\times/F}
\GL(V_a^{F^{\deg(a)}}) \cong \prod_{\bar a \in K^\times/F}
\GL_{|\gamma(\bar a)|}(q^{\deg(a)}).
\end{equation} We write
$H(\gamma) = C_G(s(\gamma))$.

The action of the unipotent part $u(\gamma)$ on $V_a$ is given by
the Jordan matrix $J(1,\gamma_0(a))$ with respect to the basis
$\BB(\gamma)_a$ of $V_a$.

\medskip

Next we define an equivalence relation on $\Gamma$ that gives rise
to the desired partition of the conjugacy classes of $G$. For
$\gamma, \delta \in \Gamma$, we write $\gamma \sim \delta$ if there
is a degree preserving bijection $\Upsilon : K^\times/F \to
K^\times/F$ such that $\gamma = \delta\Upsilon$. This defines an
equivalence relation on $\Gamma$ and for $\gamma,\delta,\Upsilon$ as
above we say $\gamma \sim \delta$ via $\Upsilon$.

For fixed $q$, the equivalence classes of $\sim$ are parameterized
by maps
\[\psi: \ZZ_{\ge 1} \to \MultiPart,
\]
written
\begin{equation} \label{e:psi(j)}
\psi(j) =
(\psi(j)_1^{b(j)_1},\psi(j)_2^{b(j)_2},\dots,\psi(j)_{m(j)}^{b(j)_{m(j)}})
\end{equation}
such that:
\begin{enumerate}
\item[(i)]
$\psi(j)$ is the empty multi-partition for all but finitely many $j
\in \ZZ_{\ge 1}$;
\item[(ii)]
$ \sum_{j \in \ZZ_{\ge 1}} j |\psi(j)| = n$; and
\item[(iii)]
$\sum_{r=1}^{m(j)} b(j)_r \le \phi(j)$ for all $j \in \ZZ_{\ge 1}$,
where $\phi$ is as in \eqref{e:Fsharp}.
\end{enumerate}
We write $\Psi$ for the set of all maps $\psi: \ZZ_{\ge 1} \to
\MultiPart$ satisfying conditions (i) and (ii) above.  For $\psi \in
\Psi$ written as in \eqref{e:psi(j)} we define
\begin{equation} \label{e:A(psi)}
A(\psi) = \{(j,r,s) \mid j \in \ZZ_{\ge 1}, r=1,\dots,m(j),
s=1,\dots,b(j)_r\}.
\end{equation}
Provided condition (iii) above holds for $\psi \in \Psi$, we may
choose $\bar a(j)_r^s \in \FF_{q^j}^\sharp/F$ for each $(j,r,s) \in
A(\psi)$ such that the $\bar a(j)_r^s$'s are all distinct. Then we
may define $\gamma \in \Gamma$, by
\begin{equation}
\label{e:gamma(a)}
\gamma(\bar a) =  \left\{ \begin{array}{cl}
            \psi(j)_r & \text{if $\bar a = \bar a(j)_r^s$, for some $(j,r,s) \in A(\psi)$}; \\
                  () & \text{otherwise}.
                \end{array}
                  \right.
\end{equation}
All possible choices for the $\bar a(j)_r^s$ gives the
$\sim$-equivalence class $\tilde \psi$ corresponding to $\psi$. If
condition (iii) does not hold for $\psi$, then, by convention, $\tilde
\psi$ is the empty set.  With this convention we can view the set
$\Psi$ as parameterizing the equivalence classes of $\sim$, and this
parametrization does not depend on $q$.

Next we count the number of elements in $\tilde \psi$ for $\psi \in
\Psi$. If we write $\psi(j)$ as in \eqref{e:psi(j)}, then,
using the description of the equivalence class $\tilde \psi$ as given
by \eqref{e:gamma(a)}, one can
see that the desired number is
\begin{equation}
\label{e:equiv} |\tilde \psi| = \prod_{j \in \ZZ_{\ge 1}}
\Delta(b(j),\phi(j)),
\end{equation}
where: $\Delta$ is defined in \eqref{e:Delta}; $b(j) =
(b(j)_1,\dots,b(j)_{m(j)}) \in \ZZ_{\ge 1}^{m(j)}$ as in
\eqref{e:psi(j)}; and $\phi(j) = |\FF_{q^j}^\sharp/F|$, see
\eqref{e:Fsharp}. Since each $\phi(j)$ is a polynomial in $q$ and
$\Delta(b(j),\phi(j))$ is a polynomial in $\phi(j)$, we see that
$|\tilde \psi|$ is a polynomial in $q$; we note, however, that in
general the coefficients of this polynomial are not integers.

\smallskip

If $\gamma \sim \delta$ (via $\Upsilon$), then we may identify the
bases $\BB(\gamma)$ and $\BB(\delta)$ of $K^n$ used to define
$x(\gamma)$ and $x(\delta)$, i.e.\ for $\bar a \in K^\times/F$, we
identify $\BB(\gamma)_{\bar a}$ with $\BB(\delta)_{\bar b}$, where
$\bar b = {\Upsilon}(\bar a)$.  Therefore, for $\psi \in \Psi$ we
may define $\BB(\psi) = \BB(\gamma)$ for some $\gamma \in \tilde
\psi$.   Suppose that $\gamma, \delta \in \tilde \psi$, then having
identified $\BB(\gamma) = \BB(\delta) = \BB(\psi)$, we have
$H(\gamma) = H(\delta)$.  Writing $H(\psi) = H(\gamma)$, we see from
\eqref{e:cent}
and the description of $\gamma \in \tilde \psi$ as in
\eqref{e:gamma(a)} that
\begin{equation} \label{e:H(psi)}
H(\psi) \cong \prod_{(j,r,s) \in A(\psi)} \GL_{|\psi(j)_r|}(q^j).
\end{equation}
We also have $u(\gamma) = u(\delta)$, so we may define $u(\psi) =
u(\gamma)$.  The conjugacy class of $u(\psi)$ in $H(\psi)$ is
parameterized by the partitions in the $\psi(j)$, i.e.\ the
conjugacy class of a unipotent element $u \in H(\psi)$ is given by
the class of the projection of $u$ into each factor
$\GL_{|\psi(j)_r|}(q^j)$, this is given by a partition of
$|\psi(j)_r|$, for $u = u(\psi)$ this is precisely the partition
$\psi(j)_r$.

For each value of $q$ such that $\tilde \psi$ is non-empty, we choose
some $\gamma = \gamma(q) \in \tilde \psi$. Then we set $x(\psi) =
x(\gamma)$, and allow this to vary as $q$ does; we note that
$x(\psi)$ depends on the choice of $\gamma$.  We write the Jordan
decomposition of $x(\psi)$ as $x(\psi) = s(\psi)u(\psi)$.  The
semisimple part $s(\psi)$ depends on the choice of $\gamma$, but
$H(\psi) = C_G(s(\psi))$ does not; $H(\psi)$ is given as in
\eqref{e:H(psi)} for all values of $q$. The parameterization of the
conjugacy class of $u(\psi) \in H(\psi)$ does not change as $q$
varies.  The discussion in this paragraph gives a convention to vary
$q$, which we use in the next section.

\section{Proof of Theorem \ref{T:main}} \label{sect:green}

For this section we fix an $n$-dimension vector $d$ and let $P =
\P_{n,d}(q)$ be the corresponding parabolic subgroup of $G =
\GL_n(q)$ as defined in the introduction.  Let $\psi \in \Psi$ and
assume $q$ is large enough so that $\tilde \psi$ is non-empty.  Let
$x = x(\psi)$, $s = s(\psi)$, $u = u(\psi)$, $\BB = \BB(\psi)$ and
$H = H(\psi) = C_G(s)$ be defined by choosing $\gamma \in \tilde
\psi$ as at the end of Section \ref{sect:conj}.

The basis $\BB = \BB(\psi)$ of $K^n$ determines an $F$-stable
maximal torus ${\bf T} = {\bf T}(\psi)$ of ${\bf G} = \GL_n(K)$
consisting of the elements of ${\bf G}$ which act diagonally on $K^n$ with
respect to $\BB$; we write $T = {\bf T}^F$.  We note that ${\bf T}$
is not split unless, $\psi(j) =()$ for all $j \ge 2$, but ${\bf T}$
is a maximally split maximal torus of ${\bf H} = C_{\bf
G}(s(\psi))$.

Suppose $x \in {}^y P$ for some $y \in G$.  Uniqueness of Jordan
decompositions implies that $s \in {}^y P$, which in turn implies
that ${}^y P \cap H$ is a parabolic subgroup of $H$.  It follows
that there exists $z \in H$ such that $T \sub {}^{zy} P$.

As $s$ is central in ${\bf H}$ and the centre of $\bf H$ is
connected, we have that $s$ is in any parabolic subgroup of $\bf H$.
In particular, this implies that $s \in Q$, for any parabolic
subgroup $Q$ of $H$, and so $x \in Q$ if and only if $u \in Q$.

We let $\CQ$ be a set of representatives of the $H$-orbits in
$\{{}^g P \mid g \in G\}$ that are of the form $H \cdot ({}^g P)$
for some ${}^g P$ with $T \sub {}^g P$; we assume that $T \sub P'$
for all $P' \in \CQ$. From the discussion in the previous two
paragraphs, we see that
\begin{equation}
\label{e:sumoverparas} f_P^G(x) = \sum_{P' \in \CQ} f_{P' \cap
H}^H(u),
\end{equation}
where the function $f_P^G$ is defined as in the introduction.  We
note that this equation does not depend on the choice of $\gamma \in
\tilde \psi$ used to define $x = x(\gamma)$.

Below we give a parameterization of the set $\CQ$. This is first
done in terms of the chosen $\gamma \in \tilde \psi$ and then we
explain how the parameterization can be described in terms of
$\psi$.  The idea is that as any $P' \in \CQ$ contains $T$;
therefore, the corresponding parabolic subgroup ${\bf P'}$ of ${\bf
G}$ (containing ${\bf T}$ and so that $P' = ({\bf P'})^F$) is the
stabilizer in ${\bf G}$ of some flag $\{0\} \sub V_1 \sub \dots \sub
V_t = K^n$ with respect to the basis $\BB = \BB(\gamma)$, i.e.\ each
$V_i$ has a basis which is a subset of $\BB$. In order for ${\bf
P'}$ to be $F$-stable we require that whenever some $v \in \BB$ is
in $V_i$ then so is $F(v)$. Further, the action of $H$ allows the
basis elements in $\BB_a$ for fixed $a \in K^\times$ to be permuted.

We let $\CC = \CC(\gamma)$ be the set of all maps
\[
c :
K^\times/F \times \{1,\dots,t\} \to \ZZ_{\ge 0}
\]
such that: $\sum_{\bar a \in K^\times/F} \deg(a) \: c(\bar a,i) =
d_i$ for each $i = 1,\dots,t$; and $c(\bar a,i) \le c(\bar a,i+1)$
and $c(\bar a,t) = |\gamma(\bar a)|$ for all $\bar a \in
K^\times/F$. Given $c \in \CC$, $a \in K^\times$ and $i \in
\{1,\dots,t\}$ we define $\BB_{a,i}$ to consist of the first $c(\bar
a,i)$ elements of $\BB_a$. We define $V_i$ to have basis $\BB_i =
\bigcup_{a \in K^\times} \BB_{a,i}$. The parabolic subgroup $Q(c)$
of $G$ is defined to be the stabilizer in $G$
of the flag $\{0\} \sub V_1 \sub
\dots \sub V_t = K^n$. We may take $\CQ = \{Q(c) \mid c \in \CC\}$
to be our set of representatives.

We write $\psi(j)$ as in \eqref{e:psi(j)} and define $A(\psi)$ as in
\eqref{e:A(psi)}.  Then $\CE = \CE(\psi)$ is defined to be the set
of all maps
\[
e : A(\psi) \times \{1,\dots,t\} \to \ZZ_{\ge 0},
\]
such that: $\sum_{(j,r,s) \in A(\psi)} j \: e(j,r,s,i) = d_i$ for
all $i = 1,\dots,t$; and $e(j,r,s,i) \le e(j,r,s,i+1)$ and
$e(j,r,s,t) = |\psi(j)_r|$ for all $(j,r,s) \in A(\psi)$.  We are
assuming that $\tilde \psi$ is non-empty, so we may fix a choice of
distinct $\bar a(j)_r^s \in \FF_{q^j}^\sharp/F$, and define $\gamma$
from $\psi$ as in \eqref{e:gamma(a)}.   For each $e \in \CE$, we
define $c = C(e) \in \CC = \CC(\gamma)$ by
\begin{equation} \label{e:c(a,i)}
c(\bar a,i) =  \left\{ \begin{array}{cl}
            e(j,r,s,i) & \text{if $\bar a = \bar a(j)_r^s$, for some $(j,r,s) \in A(\psi)$}; \\
            0 & \text{otherwise}.
                \end{array}
                  \right.
\end{equation}
The map $C : \CE \to \CC$ is a bijection.  For $e \in \CE$, we set
$Q(e) = Q(C(e))$ and note that this does not depend on the choice of
$\gamma$, i.e.\ the choice of the $\bar a(j)_r^s$.  It follows that
the set $\CE$ gives a parameterization of the set $\CQ$.

\medskip

Now by \eqref{e:sumoverparas} we get
\begin{equation}
\label{e:f_P^G} f_P^G(x(\psi)) = \sum_{e \in \CE} f_{Q(e) \cap
H}^H(u(\psi)).
\end{equation}
For values of $q$ such that $\tilde \psi$ is non-empty, each $f_{Q(e)
\cap H}^H(u(\psi))$ is a polynomial in $q$ (with integer
coefficients), by \cite[Thm.\ 3.10]{GR}. Here we use the convention
to vary $q$ as discussed at the end of Section \ref{sect:notn}.  As
the set $\CE$ does not depend on $q$, we deduce that
$f_P^G(x(\psi))$ is a polynomial in $q$.

Now by \eqref{e:kPG} we have
\[
k(P,G)  =  \sum_{\gamma \in \Gamma} f_P^G(x(\gamma)),
\]
using the parameterization of the $G$-conjugacy classes given in Section
\ref{sect:conj}.  It is implicit in \eqref{e:f_P^G} that
$f_P^G(x(\gamma)) = f_P^G(x(\psi))$ for any $\gamma \in \tilde
\psi$, so we have that
\begin{equation} \label{e:kPGpsi}
k(P,G) = \sum_{\psi \in \Psi} |\tilde \psi| \: f_P^G(x(\psi)),
\end{equation}
where, by convention we set $f_P^G(x(\psi)) = 0$ if $\tilde \psi =
\varnothing$. By \eqref{e:equiv} we have that $|\tilde \psi|$ is a
polynomial in $q$ and above we have shown that $f_P^G(x(\psi))$ is a
polynomial in $q$. Hence, $k(P,G)$ is a polynomial in $q$.

\medskip

To complete the proof of Theorem \ref{T:main}, we need to show that
the coefficients of the polynomial $k(P,G)$ are integers.  We fix a
prime $p$ and in this paragraph just consider values of $q$ that are
a power of $p$; it suffices for the proof that the coefficients of
the polynomials $k(P,G)$ are integers to just consider such $q$.
Arguing as in the introduction of \cite{Go}, we may find a family of
varieties $V_1,\dots,V_m$ defined over $\FF_p$ such that the
$P$-conjugacy classes in $G$ correspond to the $\FF_q$-rational
points of the $V_i$.  More precisely, using Rosenlicht's theorem
(see \cite{Ros}), we may find a $\bf P$-stable open subvariety $U_1$
of $\bf G$ defined over $\FF_p$ and an orbit space $V_1$ for the
action of $\bf P$ on $U_1$. This means that the points of $V_1$
(over $K$) correspond to the ${\bf P}$-conjugacy classes in $U_1$.
Now using the fact that $C_{\bf P}(x)$ is connected for any $x \in
{\bf G}$, we see that the $\FF_q$-rational points of $V_1$
correspond to the conjugacy classes of $P$ in the set of
$\FF_q$-rational points of $U_1$, this follows from \cite[Prop.\
3.21]{DM}.  Now we may apply Rosenlicht's theorem to the action of
$\bf P$ on ${\bf G} \setminus U_1$ to find $U_2$ and $V_2$ in
analogy to $U_1$ and $V_1$. Continuing in this way, we obtain the
varieties $V_1,\dots,V_m$ whose $\FF_q$-rational points correspond
to the $P$-conjugacy classes in $G$. Given this parameterization of
the $P$-conjugacy classes in $G$, one can apply some standard
arguments, using the Grothendieck trace formula (see \cite[Thm.\
10.4]{DM}), to prove that the coefficients of the polynomial
$k(P,G)$ are integers, see for example \cite[Prop.\ 6.1]{Re}.

We note that the polynomial summands $|\tilde \psi| \:
f_P^G(x(\psi))$ in the expression for $k(P,G)$ given in
\eqref{e:kPGpsi} do not have integer coefficients in general; this
can already be seen for $G = \GL_2(q)$ in the examples below.

\medskip

We conclude our discussion with some examples which demonstrate that
it is possible to explicitly calculate the polynomials $k(P,G)$.  We
observe that in the examples below, $k(P,G)$ is divisible by $q-1$.
One can see that this has to be the case by checking that $q-1$
divides the polynomial $|\tilde \psi|$ for all $\psi$.

\begin{exmps}
(i).  We begin by explicitly calculating $k(B,G)$ and $k(G)=k(G,G)$
for $G = \GL_2(q)$.  The possible values of $\psi$ and all the
information needed to calculate $k(B,G)$ and $k(G)$ is given in the
table below. It is straightforward to calculate all of the
information in this table by hand.

\begin{table}[h!tb]
\renewcommand{\arraystretch}{1.5}
\begin{tabular}{|c|c|l|c|c|}
\hline $\psi(1)$ & $\psi(2)$ & \qquad\quad\ $x(\psi)$ & $|\tilde \psi|$ & $f_B^G(x(\psi))$ \\
\hline $((1^2))$ & $()$ & $\left(\begin{array}{cc} a & 0 \\ 0 & a
\end{array} \right)$, $a \in \FF_q^\times$ & $q-1$ & $q+1$\\
$((2))$ & $()$ & $\left(\begin{array}{cc} a & 1 \\ 0 & a
\end{array}\right)$, $a \in \FF_q^\times$ & $q-1$ & $1$\\
$((1)^2)$ & $()$ & $\left(\begin{array}{cc} a & 0 \\ 0 & b
\end{array}\right)$, $a \ne b \in \FF_q^\times$ & $\frac{(q-1)(q-2)}{2}$ & $2$
\\
$()$ & $((1))$ & $\left(\begin{array}{cc} a & 0 \\ 0 & a^q
\end{array}\right)$, $a \in \FF_{q^2}^\sharp$ & $\frac{q^2 - q}{2}$ &
$0$
\\ \hline
\end{tabular}
\medskip
\label{Tab:GL2}
\end{table}

Now using \eqref{e:kPGpsi} we can calculate:
\[
k(B,G) = (q-1)(q+1) + (q-1)1 + \frac{(q-1)(q-2)}{2}2 = 2q(q-1).
\]
Of course, we have $f_G^G(x(\psi))=1$ for all $\psi$, so we obtain:
\[
k(G) = (q-1) + (q-1) + \frac{(q-1)(q-2)}{2} + \frac{q^2 - q}{2} =
(q-1)(q+1).
\]

(ii).  For $n \ge 3$ (not too large), it is straightforward to
calculate $k(B,G)$ using the values of the functions $f_B^G(u)$ for
unipotent $u$.   It is possible to obtain these values using the
{\tt chevie} package in GAP3 (\cite{gap}) along with some code
provided by M.~Geck, and the formula for $f_B^G(u)$ given in
\cite[Lem.\ 3.2]{GR}. The size of $\Psi$ gets large quickly as $n$
increases, so we have only calculated the values of $k(B,G)$ for $n
\le 4$.  We do not include the details of these calculations here,
as that would take a lot of space.  For $n=3$, we get
\[
k(B,G) = (q-1)(q^3+6q^2-q-3)
\]
and for $n = 4$ we obtain
\[
k(B,G) = (q-1)(q^6+3q^5+9q^4+19q^3-9q^2-18q+5).
\]

(iii).  We finish by giving an example of how to calculate a
particular value of $f_P^G(x(\psi))$.  We consider the case $G =
\GL_9(q)$, $P = \P_{9,d}(q)$, where $d$ is the $9$-dimension vector
$(4,7,9)$, and $\psi$ is given by
$$
\psi(1) = ((2)), \, \psi(2) = ((1^2)), \, \psi(3) = ((1)); \, \text{
and } \psi(j) = () \text{ for } j \ge 4.
$$

We write $x = x(\psi)$ with Jordan decomposition $x = su$, and we
write $H = C_G(s)$.  We have the direct product decomposition $H =
\GL_2(q) \times \GL_2(q^2) \times \GL_1(q^3) = H_1 \times H_2 \times
H_3$ say. We write $x_i$ for the projection of $x$ into $H_i$, for
each $i$.  We note that $x_1$ is a product of a central element and
a regular unipotent element in $H_1$, $x_2$ is central in $H_2$ and
$x_3$ is central in $H_3$. Given a parabolic subgroup $Q$ of $H$
containing $s$, we write $Q_i = Q \cap H_i$ for each $i$, and note
that
\begin{equation} \label{e:prod}
f_Q^H(x) = f_{Q_1}^{H_1}(x_1)f_{Q_2}^{H_2}(x_2)f_{Q_3}^{H_3}(x_3).
\end{equation}
Using \eqref{e:equiv} we can calculate
\[
|\tilde \psi| = (q-1)\frac{q^2-q}{2}\frac{q^3-q}{3}.
\]
We have $A(\psi) = \{(1,1,1),(2,1,1),(3,1,1)\}$.  There are three
elements $e \in \CE(\psi)$, they are shown in the following three
matrices: the value of $e(j,1,1,i)$ being given by the entry in the
$j$th row and $i$th column:
\[
\left(\begin{array}{ccc}
1 & 2 & 2\\
0 & 1 & 2\\
1 & 1 & 1 \\
\end{array}\right)
\hspace{10mm} \left(\begin{array}{ccc}
0 & 0 & 2\\
2 & 2 & 2\\
0 & 1 & 1 \\
\end{array}\right)
\hspace{10mm} \left(\begin{array}{ccc}
2 & 2 & 2\\
1 & 1 & 2\\
0 & 1 & 1 \\
\end{array}\right)
\]

Next we use \eqref{e:prod} to work out the value of $f_{Q(e)}^H(x(\psi))$
for each of the three possible values of $e$.  In the first
case we have that $Q_1$ is a Borel subgroup of $H_1$, so that
$f_{Q_1}^{H_1}(x_1) = 1$;  $Q_2$ is a Borel subgroup of $H_2$, so
that $f_{Q_2}^{H_2}(x_2) = q^2+1$; and $Q_3$ is (necessarily) all of
$H_3$, so we get $f_{Q_3}^{H_3}(x_3) = 1$. We can work out the value
of $f_{Q(e)}^H(x)$ for the other two possible values of $e$
similarly, and then we can use \eqref{e:f_P^G} to calculate
\[
f_P^G(x) = (q^2+1) + 1 + (q^2+1) = 2q^2+3.
\]
\end{exmps}

%%%%%%%%%%%%%%%%%%%%%%%%%%%%%%%%%%%%%%%%%%%%%%%%%%%%%%%%%%%%%%%%%%%%%
%%%%%%%%%%%%% Acknowledgments
%%%%%%%%%%%%%%%%%%%%%%%%%%%%%%%%%%%%%%%%%%%%%%%%%%%%%%%%%%%%%%%%%%%%%%

\bigskip

{\bf Acknowledgments}: This research was funded in part by EPSRC
grant EP/D502381/1.  The first author would like to thank New
College, Oxford for financial support whilst part of this research
was carried out.  We thank the referees for useful comments and
suggestions.

%%%%%%%%%%%%%%%%%%%%%%%%%%%%%%%%%%%%%%%%%%%%%%%%%%%%%%%%%%%%%%%%%%%%%%
%%%%%%%%%%%%% bibliography
%%%%%%%%%%%%%%%%%%%%%%%%%%%%%%%%%%%%%%%%%%%%%%%%%%%%%%%%%%%%%%%%%%%%%%

\end{document}